\documentclass[10pt,amstex]{article}
\usepackage[usenames]{color}

\usepackage{amscd}
\usepackage{amsmath}
\usepackage{amssymb}
\usepackage{graphicx}

\newtheorem {theo} {\bf Theorem} [section]
\newtheorem {prop} [theo] {\bf Proposition}

\newtheorem {defn} [theo] {\bf Definition}

\newtheorem {rem} [theo] {\bf Remark}

\newcommand{\qed}{\nopagebreak\hfill{\vrule width6pt height6pt depth0pt}}
\newcommand{\be}{\begin{eqnarray}}
\newcommand{\ee}{\end{eqnarray}}
\newcommand{\benn}{\begin{eqnarray*}}
\newcommand{\eenn}{\end{eqnarray*}}
\newcommand{\bse}{\begin{equation}}
\newcommand{\ese}{\end{equation}}
\newcommand{\bsenn}{\begin{displaymath}}
\newcommand{\esenn}{\end{displaymath}}
\newcommand{\logand}{\;\;{\rm and }\;\;}

\newcommand{\where}{\;\;{\rm where }\;\;}

\newcommand{\C}{\mathbb{C}}
\newcommand{\N}{\mathbb{N}}

\newcommand{\R}{\mathbb{R}}

\usepackage{graphicx}

\begin{document}

\title{Transients in the Synchronization of Oscillator Arrays}
\author{ Carlos E. Cantos\thanks{Fariborz Maseeh Dept. of Math. and Stat., Portland State Univ., Portland, OR, USA.}, J. J. P. Veerman\thanks{Fariborz Maseeh Dept. of Math. and Stat., Portland State Univ., Portland, OR, USA, and also CCQCN, Dept of Physics, University of Crete, 71003 Heraklion, Greece; e-mail: veerman@pdx.edu}\\
}\maketitle

\vskip .0in

 %\normalsize        %\mysetfontsize5

\noindent
\section*{Abstract}
The purpose of this note is threefold. First we state a few conjectures that allow us to
rigorously derive a theory which is asymptotic in $N$ (the number of agents) that describes
transients in large arrays of (identical) linear damped harmonic oscillators in $\R$ with
completely decentralized nearest neighbor interaction.
We then use the theory to establish that in a certain range of the parameters transients
grow \emph{linearly} in the number of agents (and faster outside that range). Finally, in the regime
where this linear growth occurs we give the constant of proportionality as a function
of the signal velocities (see \cite{CVH}) in each of the two directions. As corollaries
we show that symmetric interactions are far from optimal and that all these results
independent of (reasonable) boundary conditions.

%\vskip 0.4in\noindent
\begin{centering}
\section{Introduction}
\label{chap:zero}
\end{centering}
\setcounter{figure}{0} \setcounter{equation}{0}

In this paper we formulate a theory for transients of certain large arrays of linear damped harmonic
oscillators. For our purposes transients are solutions of a high-dimensional system $\dot \zeta= f(\zeta)$
that converge to a stable equilibrium. We investigate transients of the linearized equations, i.e. $\dot \zeta=M\zeta$ where $M$ is a square matrix.
It is well-known (see \cite{Tr} and \cite{TTRD}) that transients of high dimensional systems
cannot successfully be analyzed by looking at the spectrum of the linear operators unless the
operator is \emph{normal}.
An $N\times N$ square matrix $M$ is \emph{normal} if it commutes with its (conjugate) transpose, or equivalently
(by the spectral theorem), if its eigenvectors form an orthonormal basis of $\C^N$.
If the operator $M$ is not normal then only in the limit as $t\rightarrow \infty$ is the
spectrum decisive.
The spectrum of the operator gives therefore very little information about high dimensional transients.

In this work we consider arrays of (linear damped harmonic) oscillators on the line interacting with their nearest neighbors.
We label the agents from 0 to $N$ (from right to left) and describe the interaction as, for all $k\in\{0,\cdots N\}$:
\bse
\ddot x_k = f(x_k-x_{k-1}, x_k-x_{k+1},\dot x_k - \dot x_{k-1},\dot x_k-\dot x_{k+1})
\label{eq:nonlinear}
\ese
The orbit, $x_0(t)$,  of \emph{the leader} is prescribed as $x_0(t)=\max\{0,v_0t\}$ (to the right).
We formulate a theory that quantitatively describes the response of the flock if $N$ is large.
Before we do so, it will serve us well to keep one of the main applications of our theory in mind.
Suppose we have a long sequence of cars equipped with automatic pilots waiting for
a traffic light. At $t=0$ the light turns green and the lead car (or \emph{leader}) acquires constant speed $v_0>0$.
Each car is equipped with a sensor that perceives relative position and velocity
of each of their neighbors (back and front). This input is then used as feedback
for the acceleration of each of the cars. The purpose is that the cars follow each
other as well as possible, until everyone moves behind the leader with the same speed.
It turns out that if done correctly, the acceleration of the leader causes a perturbation
throughout the system (the \emph{transient}), until the system, assuming it is asymptotically stable, settles down again. The question is: What is the algorithm for the feed-back that minimizes the transients? While a lot of work has been done to find interesting systems of differential
equations that can simulate flock behavior --- many of these consist of
nonlinear equations (see \cite{LF}, \cite{OS}) --- our emphasis is more related to performance.
We aim to characterize the transient qualitatively and quantitatively when the number of
cars is large.

We formulate a theory (for the indicated initial conditions) that is comprehensive in that it
covers all parameter values of the models we consider (including the different boundary conditions),
explicit in the sense that it gives explicit answers for what the transients are,
and simple in the sense that it generalizes to more complex situations.
The cases covered by the theory developed
in this paper cover all possible linearizations of Equation \ref{eq:nonlinear}
(and that includes non-symmetric interactions and different boundary conditions).
To the best of our knowledge there is no literature that considers (let alone analyzes)
the problem in this generality.
The price we pay for such a general theory is that we cannot prove all our statements.
We need some conjectures, and these will be spelled out in Section \ref{chap:compare}.

In Section \ref{chap:circulant} we set the definitions and sketch some of our previous results.
Those results do not depend on conjectures, but are more restrictive and the methods are hard to generalize.
We also introduce the notion of \emph{flock stability}.
Flock stability means that (time-)responses to initial perturbations grow less than
exponential in the number $N$ of agents as $N$ tends to infinity.
This notion is independent of that of \emph{asymptotic stability} which deals with
the growth of the response of a single system (with $N$ fixed) the response as $t\rightarrow\infty$.

Section \ref{chap:compare} contains our main results.
In Definition \ref{def:principles} we set out a number
of principles (conjectures) that allow us to rigorously derive a theory for these systems.
Using this theory we then find the parameter range for which parameters the transients are
well-behaved (Proposition \ref{prop:stable}). Within that range, solutions of our system
are traveling waves with high frequency attenuation, and we give a precise analytical description
of the transients as functions of the parameters (Theorem \ref{theo:quantities}). The theory
is asymptotic in the number of agents.

Since our theory rests on conjectures we need to verify its conclusions.
In Section \ref{chap:numerical} we compare our results with measurements done in 360
simulations. The results are in excellent agreement. We close with some final remarks
in Section \ref{chap:five}.

%\vskip 0.4in
\begin{centering}\section{Nearest Neighbor Flocks in $\R$}
 \label{chap:circulant}\end{centering}
\setcounter{figure}{0} \setcounter{equation}{0}

Ever since the inception (\cite{HM1}, \cite{HM2}) of the subject it has been
a challenge to mathematically express the notion \emph{decentralized} (\cite{Chu}).
Here is how we interpret it (following \cite{flocks2} and \cite{position}): In a
decentralized flock the only information an agent receives are the position and velocity
relative to it of nearby agents, i.e. the acceleration of the $k$th agent is a function
of the differences $x_k-x_i$ and $\dot x_k - \dot x_i$ where $i$ runs over the neighbors
of $k$ and agents do not receive information from an outside source.
The only exception to this is the \emph{leader}, whose orbit is prescribed.
This motivates the following definition.

\vskip .0in\noindent
\begin{defn}
We consider the system given by Equation \ref{eq:nonlinear} for all $k\in\{0,\cdots N\}$,
where for $k\in\{1,\cdots N-1\}$ the $f_k$ are equal to $f:\R^4\rightarrow \R$, $f_0=0$, and
$f_N$ is different from the others due to boundary conditions. In addition we require that
there is $\Delta>0$ so that $f_k(-\Delta, \Delta, 0,0)=0$.
\label{defn: nearest neighbor}
\end{defn}

%\noindent
For all real constants $x_0$ and $v_0$ and for positive $\Delta$
\bsenn
\forall\;k\in\{0,\cdots N\} \, : \quad  x_k = v_0 t+x_0-\Delta k
\esenn
are the \emph{coherent motions} of the flock as a whole with velocity $v_0$ and inter-agent distance $\Delta$.
It is easy to check that these are solutions of the Equation \ref{eq:nonlinear} if and only if
for all $k$ we have $f_k(-\Delta, \Delta, 0,0)=0$. Hence the additional requirement in the definition.

If we substitute
\bsenn
x_k=v_0t+x_0-\Delta k + \varepsilon z_k
\esenn
into Equation \ref{eq:nonlinear} and work out the first order in $\varepsilon$ and use that $f_k(-\Delta, \Delta, 0,0)=0$. we then obtain the following Proposition.

\vskip .0in
\begin{prop} The linearization of the system given in Definition \ref{defn: nearest neighbor} around the coherent solutions is, for all $k\in\{1,\cdots N-1\}$:
\be
\ddot z_k & = & g_x \left[\rho_{x,-1} z_{k-1} + \rho_{x,0}z_k+\rho_{x,1} z_{k+1}\right] + g_v \left[\rho_{v,-1}\dot z_{k-1} + \rho_{v,0}\dot z_k +\rho_{v,1} \dot z_{k+1}\right]
\label{eqn:oscillator array2}
\ee
Here $g_x$, $g_v$, $\rho_{x,i}$ and $\rho_{v,i}$ are constants, $\sum_j\,\rho_{x,j}=0$, and $\sum_j\,\rho_{v,j}=0$.
At $t=0$ the system is at rest and after that the leader moves with constant velocity ($z_0(t)=v_0t$).
The equation for the last agent is discussed in Definition \ref{defn:boundary condns}.
\label{prop:linearization}
\end{prop}

\vskip .0in
This leads us to the consideration of boundary conditions.
The orbit of the leader is prescribed. But for agent $N$ (the last one) we
have to make a choice. Some reasonable choices (but not the only ones) are enumerated in the
Definition below. If agents have more than 2 neighbors the number of `reasonable' boundary
conditions becomes greater.

\begin{defn} {\rm (Boundary Conditions)} Let $S_N$ be the linearized system in Proposition
\ref{prop:linearization}. We specify the orbits of the first and last agents. The first agent's (the leader's) motion is described by:
\bsenn
z_0(t)=\left\{ \begin{array}{cc} 0 & \rm{\;if\;} t<0\\ tv_0 & \rm{\;if\;} t\geq 0 \end{array} \right.
\esenn
The boundary conditions are expressed in the equation for $\ddot z_N$ and
are characterized by the parameters
$\beta_x$ and $\beta_v$:
\bsenn
\ddot z_N = g_x \beta_x \left[- z_{N-1} + z_N\right]
                    +g_v \beta_v \left[- \dot z_{N-1} + \dot z_N\right]
\esenn
\begin{enumerate}
\item Variable mass boundary conditions: $\beta_x=-\rho_{x,-1}$ and $\beta_v=-\rho_{v,-1}$
\item Regular boundary conditions: $\beta_x=\beta_v=1$
\end{enumerate}
\label{defn:boundary condns}
\end{defn}

\vskip .0in
The first set of boundary conditions arises from simply leaving out the
dependence on the relative velocity and position of the rear-neighbor.
When $\rho_{x,1}=\rho_{x,-1}$ and $\rho_{v,1}=\rho_{v,-1}$, these
boundary conditions give rise to \emph{symmetric} Laplacian matrices $L_x$ and $L_v$.
For this reason that boundary condition is most often used (e.g. see \cite{YDR},
\cite{PB}).
Notice though that physically this is akin to changing the mass of the last car.
In that sense the second
set of boundary conditions (used by other authors, e.g. \cite{HM1}) is more realistic.

\vskip .0in\noindent
We now give a definition of flock stability that is sufficient for our considerations.
(More details can be found in \cite{flocks6, flocks7, flocks67}.) The definition of asymptotically
stable is standard (see \cite{Ro}).

\begin{defn} {\rm (Flock stability)} The collection of systems $S_N$ of Definition
\ref{defn:boundary condns} is called flock stable if it is asymptotically stable and
if $\max_{t\in\R}|z_N(t)|$ grows sub-exponentially in $N$.
\label{defn:flock-stable}
\end{defn}

\vskip .0in\noindent
Earlier characterizations of flocks concentrated on \emph{string stability}, which is the variation
of the \emph{relative distances between neighbors} (e.g. \cite{Chu}, \cite{Pe}, \cite{Co},
\cite{SHCI}, \cite{PSH}, \cite{SH1}, \cite{SH2}, \cite{MB}).
This is hard to generalize for more complicated flocks in dimension 2 or higher.
Also most of these papers, and many others (\cite{MB}, \cite{YDR}), in fact consider the size of
frequency response as a measure of stability. While this is mathematically equivalent
to the time-domain response, it is often difficult or impossible to calculate the time-response
from the frequency response (see for instance \cite{flocks6}, \cite{flocks7}, \cite{HMHS},
\cite{MHHS}). Since it is our interest to calculate precise shapes and sizes of transients
in the time-domain, we use a very different approach.

\begin{figure}[pbth]
\center
\includegraphics[width=.4\linewidth]{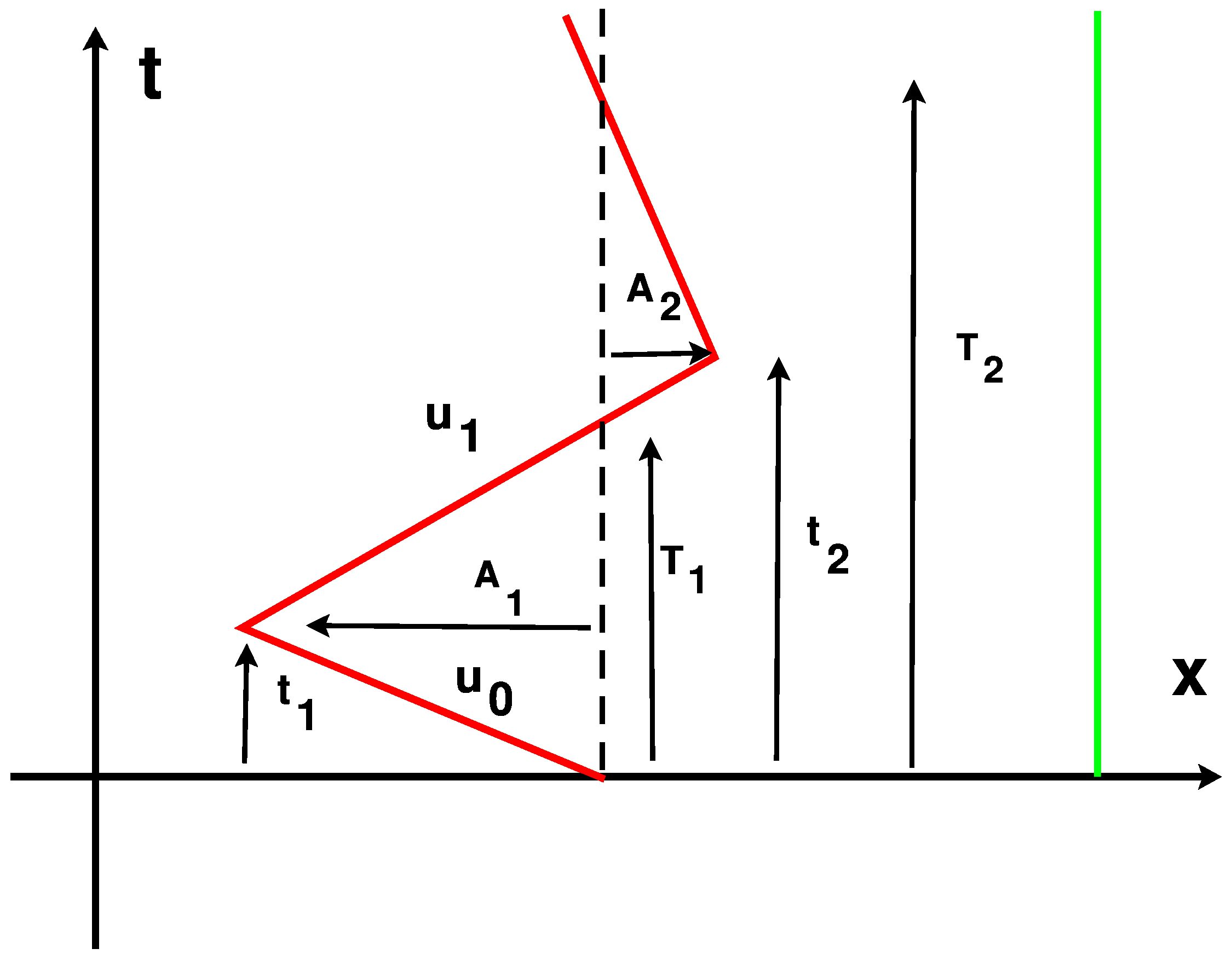}
\caption{ \emph{ Schematic representation of orbit of agent with largest fluctuation relative
to the leader. Horizontal is position, vertical is time. The vertical (green) segment on the far right
represents the leader's orbit. The dashed line represents the equilibrium position of the trailing agent. The oscillating (red) line represents the actual orbit. The response times $t_{i}$,
the times $T_k$ at which the orbit crosses its equilibrium
position, and the amplitudes $A_k$ are indicated. }}
\label{fig:schematic}
\end{figure}

%\vskip .0in
We now summarize what is known about the transients of the system in Definition
\ref{defn: nearest neighbor}. In Figure \ref{fig:schematic}
a simulation of the response when the leader leaves at $t=0$ with unit velocity is shown.
Each represents the orbit of an agent \emph{relative to the leader} (so that the leader
appears to stand still).
In what follows we will characterize the orbit of the last agent, $z_N(t)$, in terms of its
local extrema $A_k$, the period $T$, which is the (average) time elapsed between $T_k$ and $T_{k+2}$,
and the attenuation $\alpha$, which is the (average of the) ratio $A_{k+2}/A_k$ (see Figure \ref{fig:schematic}).

\begin{theo} i) (\cite{flocks5}) The system $S_N$ (see Definition \ref{defn:boundary condns})
with  $\rho_{v,1}=\rho_{x,1}=-1/2$ and regular boundary conditions, is flock stable and
\bsenn
\alpha = 1 \quad\logand \quad  A_1= -\sqrt{\dfrac{|g_x|}{2}} \quad \logand \quad
T_2=\dfrac{4\sqrt2 N}{\sqrt{|g_x|}}
\esenn
All expressions are asymptotic in $N$.\\
ii) (\cite{flocks6}, \cite{flocks7}, \cite{flocks67})
The system $S_N$ (see Definition \ref{defn:boundary condns})
with  $\rho_{v,1}=\rho_{x,1}=r$ and regular boundary conditions,
is asymptotically stable and flock unstable for all $r\in(-1,0)\backslash \{-1/2\}$.
\label{theo:flock5}
\end{theo}

\vskip .0in \noindent
{\bf Remark:} We note that the case $r\in(-1/2,0)$ is particularly interesting since the real parts
of the non-zero eigenvalues are bounded by
$\max\{-\dfrac{g_x}{g_v}, g_v(1-2\sqrt{|r|(1+r)})\}$
while when $r=-1/2$ the real parts tend to 0 as $N\rightarrow \infty$.
Nonetheless the only flock stable case is the latter.

%\vskip .0in\noindent
Finally we need the nearest neighbor model $S_N^*$ with periodic boundary conditions (no leader) studied in \cite{CVH}.

%\vskip .0in
\begin{defn} The standard system $S_N^*$ is given as follows. For all $k\in\{1,\cdots N\}$,
\benn
\ddot z_k & = & g_x \left[\rho_{x,-1} z_{k-1} + \rho_{x,0}z_k+\rho_{x,1} z_{k+1}\right] + g_v \left[\rho_{v,-1}\dot z_{k-1} + \rho_{v,0}\dot z_k +\rho_{v,1} \dot z_{k+1}\right]
\eenn
where
\bsenn
\rho_{x,j}=\rho_{x,j+N} \quad \logand \quad \rho_{v,j}=\rho_{v,j+N}\quad \logand \quad
z_{j}=z_{j+N}
\esenn
Here $g_x$, $g_v$, $\rho_{x,i}$ and $\rho_{v,i}$ are constants, and $\sum_j\,\rho_{x,j}=0$, and $\sum_j\,\rho_{v,j}=0$.
Further, let the initial conditions $z_k(0)$ and $\dot z_k(0)$ be such that the Fourier series of $z_k(0)$
has coefficients $a_m$ that decay (at least) quadratically in $m$.
\label{defn:normalized system}
\end{defn}

\begin{theo} (See \cite{CVH}) Fix $0<\alpha<1/2<\beta<2/3$ and $K>1$. Denote the solution 
of $S_N^*$ as $z_j(t)$.  For all $\varepsilon>0$ there is an $N_0 \in \N$ such 
that for any $N>N_0$ there exist waves $f_+$ and $f_-$ such that:
\bse
| z_j(t) - f_-(j-c_- t) - f_+(j-c_+ t) | < \varepsilon ,
\label{eq:traveling-wave}
\ese
for all $t \in \left[ \frac{N}{|c_-|}, K \frac{N}{|c_-|} \right] \cap \left[ \frac{N}{c_+}, K \frac{N}{c_+} \right]$.
Furthermore the \emph{signal velocities} $c_\pm$ (in terms of number of agents per unit time) are
given by:
\bse
c_{\pm} = -\dfrac{g_v(1+2\rho_{v,1})}{2}\pm \sqrt{\dfrac{g_v^2(1+2\rho_{v,1})^2}{4}-\dfrac{g_x}{2}},
\label{eq:signal-velocities}
\ese
where $c_+>0$ and $c_-<0$.
\label{theo:CV1}
\end{theo}

%\vskip .4in
\begin{centering}\section{Transients in Nearest Neighbor Systems}
\label{chap:compare}
\end{centering}
\setcounter{figure}{0} \setcounter{equation}{0}

Our aim in this section is first of all to set out some basic
principles that allow us to rigorously analyze the transients of the system. And secondly
to carry out that analysis and give a complete characterization of the transients.

\vskip .0in\noindent
\begin{defn} Let $P_N$ be the system given by
Definition \ref{defn:boundary condns} with the of exception $z_0(t)$, which is given by
$z_0(t)=p(t)v_0$. Here $p$ satisfies $\int p(s)\,ds=1$ and has compact support.
$v_0$ is the velocity of the leader (see Proposition \ref{prop:linearization}).
\label{def:second deriv}
\end{defn}

\vskip .0in\noindent
{\bf Remark:} $P_N$ is obtained by differentiating $S_N$ twice. Therefore solutions of
of $S_N$ can be obtained by integrating solutions of $P_N$ twice.

\vskip .0in\noindent
We will use the following principles to guide our solution of the system $P_N$ in the order that they will be used. We will use the word \emph{unit pulse} for a function $p(t)$ satisfying the requirements
given in Definition \ref{def:second deriv}.

\vskip .0in\noindent
\begin{defn}, {\bf Principles for solving $P_N$:} Fix all parameters and a pulse $p$
with support in an interval $[-\varepsilon,\varepsilon]$. For $N$ large enough, the following holds:
\begin{enumerate}
\item If $S_N^*$ is asymptotically unstable then $S_N$ is asymptotically unstable or
flock unstable.
\item The evolution of $P_N$ is the same as that of $S_N^*$ except at the boundaries.
\item If $P_N$ is asymptotically stable and flock stable then the left boundary condition of the solution of $P_N$ is given by $\partial_jz_{j}(t)|_{j=N}=0$.
\item The right boundary condition of the solution of $P_N$ is given by $z_0(t)=p(t)$.
\item If $P_N$ is asymptotically stable and flock stable then a unit pulse $p$ propagating
through the system remains a unit pulse until boundary effects take place.
\end{enumerate}
\label{def:principles}
\end{defn}

\vskip .0in\noindent
{\bf Remarks:} It is worthwhile to comment on the mathematical status of these principles.
Principle 4 is true, because is part of the definition of $P_N$.
We have not been able to encounter a formal proof of Principle 2. However, it seems
reasonable that a traveling signal or pulse does not ``feel'' the boundary if it is far away
from it. This principle is widely used in many areas of science (eg phonons in condensed
matter physics) and we feel justified in following the literature.
This leaves Principles 1, 3, and 5 as our real Conjectures.
The final justification of these assumptions lies in the agreement between their implied
results and numerical simulations (Section \ref{chap:numerical}).
However below we give some additional reasons why these particular conjectures make sense.

Principle 1 is based on the following intuition. It is easy to see that the $2N$-dimensional space
of solutions of the linear operator of $S_N^*$ has $N$ two-dimensional ($\mathbb{C}^2$)
eigenspaces that are orthogonal each other. The two eigenvectors within each plane are typically
not orthogonal. In changing the boundary condition from those in $S_N^*$ to those of $S_N$,
these planes lose orthogonality. Thus one expects the dynamics typical for high dimensional
non-normal systems to be more pronounced in $S_N$ than in $S_N^*$. The salient aspect of this kind
of dynamics is the fact that asymptotically stable systems nonetheless
have exponentially large transients (see \cite{TTRD} and \cite{Tr}). Since these systems
are in some sense close to another, one expects the initial behavior of both systems
to be close for some time. Thus initial growth of transients in one will indicate
growth of transients in the other: this growth may be due to asymptotic instability
or to flock instability.

Principle 3 is related to the fact that agent $N$ in the boundary
has no restrictions of movement (as opposed to agent 0 which is held fixed, see Principle 4).
Thus Principle 3 expresses what is known as a \emph{free boundary condition} in the analysis
of the wave equation in one dimension (see \cite{CH}, Appendix 2 to Chapter 5). Note that
it renders the problem \emph{independent of the details of the boundary condition} in Definition
\ref{defn:boundary condns}. (See \cite{WS} for a study of the dependence on boundary conditions.)
This free boundary condition has been used when replacing a similar system with a PDE (\cite{BMH}).
Very recently this boundary was independently proposed by \cite{MHHS} for a discrete system.

If we combine Principle 2 with Theorem \ref{theo:CV1} we see that Principle 5 holds for
low-frequency pulses. If $p(t)$ does contain high frequencies, we do not have a formal proof.
However it appears that in this case the high-frequency damping does not affect $\int z_j(t)dt$
as $j$ increases. Writing the
above solution of the wave equation as a sum of eigensolutions, it appears that the integral
only depends on coefficient of the eigensolution associated to the eigenvalue zero
(ie: it doesn't decay).

\vskip .0in\noindent
\begin{prop} A necessary condition for the system $S_N$ to be asymptotically stable and flock
stable is
\bsenn
\rho_{x,-1}=\rho_{x,1}\quad \logand \quad g_x\rho_{x,0}<0 \quad\logand \quad g_v\rho_{v,0}<0
\esenn
\label{prop:stable}
\end{prop}

\vskip -.3in\noindent
{\bf Proof:} This is a direct consequence of principle 1 and the fact that $S_N^*$
is asymptotically stable if and only if these conditions hold (see \cite{CVH}).
$\blacksquare$

\vskip .0in\noindent
{\bf Remark:} From now on we will restrict ourselves to systems that satisfy the conditions
of this proposition. To simplify notation (and without loss of generality) we will re-scale $g_x$ and
$g_v$ so that $\rho_{x,0}=\rho_{v,0}=1$. We will write:
\be
\rho_{x,0}=\rho_{v,0}=1 \quad \logand \quad \rho_{x,-1}=\rho_{x,1}=-1/2
\logand \quad  g_x<0 \quad\logand \quad g_v<0
\label{eq:new-conds}
\ee

\vskip .0in\noindent
\begin{prop} Let $P_N$ be as in Definition \ref{def:second deriv} and Equation \ref{eq:new-conds}. 
Fix the impulse $p(t)$, and $K>1$. As $N$ tends to infinity, the orbit of the last agent, in the time interval 
$t \in \left[ \frac{N}{|c_-|}, K \frac{N}{|c_-|} \right] \cap \left[ \frac{N}{c_+}, K \frac{N}{c_+} \right]$ 
will tend to:
\bsenn
z_N(t)= \dfrac{c_+-c_-}{c_+} \sum_{k=0}^{\infty} \left(\dfrac{c_-}{c_+}\right)^k
f_k\left(t-\dfrac{N}{c_+}-\left(\dfrac{1}{c_+}-\dfrac{1}{c_-}\right)kN\right) v_0
\esenn
Here the $f_k$ are unit pulses.
\label{prop:solution}
\end{prop}

\vskip .0in\noindent
{\bf Proof:} By Principle 2 the system evolves as $S_N^*$.
By Theorem \ref{theo:CV1} the solution $z_j(t)$ can be expressed as a sum of traveling waves.

It remains to analyze the effect of the boundaries at $j=N$ and $j=0$ on the solution.
For this we give a reasoning analogous to the analysis in \cite{CH}, Appendix 2 to Chapter 5.
There are two substantial modifications. The spatial variable is discrete and the signal
velocities depend on the direction. Because the system is linear we may use 1 instead of $v_0$ and multiply
the solution we then obtain by $v_0$.

We start by noting that Principle 3 gives:
\bsenn
-\dfrac{1}{c_+}\psi'\left(t-\dfrac{N}{c_+}\right)-
\dfrac{1}{c_-}\phi'\left(t-\dfrac{N}{c_-}\right)=0
\esenn
Assume that $\phi$ and $\psi$ are continuous and integrate with respect to $t$ to get:
\bse
{c_-}\psi\left(t-\dfrac{N}{c_+}\right)+{c_+}\phi\left(t-\dfrac{N}{c_-}\right)=0
\label{eq:1}
\ese
Substitute $s_-=t-N/c_-$ into Equation \ref{eq:1} to get
\bse
\phi\left(s_-\right)=-\dfrac{c_-}{c_+}\,\psi\left(s_--\left(\dfrac{1}{c_+}-\dfrac{1}{c_-}\right)N\right)
\label{eq:3}
\ese
Principle 4 and Theorem \ref{theo:CV1} (with $j=0$) give $\psi(s_-)+\phi(s_-)=p(s_-)$.
Substitute Equation \ref{eq:3} into this and we get a recursion
\bsenn
\psi(s_-)=p(s_-)+\dfrac{c_-}{c_+}\, \psi\left(s_--\left(\dfrac{1}{c_+}-\dfrac{1}{c_-}\right)N\right)
\esenn
And that implies:
\bsenn
\psi(s_-)=p(s_-)+\sum_{k=1}^{\infty}\,\left(\dfrac{c_-}{c_+}\right)^k\,
 p\left(s_--\left(\dfrac{1}{c_+}-\dfrac{1}{c_-}\right)kN\right)
\esenn
We note that for finite $s_-$, this is a finite sum because $c_-<0$ and therefore $\left(\dfrac{1}{c_+}-\dfrac{1}{c_-}\right)$ is positive, so that for $k$ sufficiently large
$p\left(s_--\left(\dfrac{1}{c_+}-\dfrac{1}{c_-}\right)kN\right)=0$.

On the other hand by using $s_+=s_--\left(\dfrac{1}{c_+}-\dfrac{1}{c_-}\right)N$, we see that Equation \ref{eq:3} gives:
\bse
\psi\left(s_+\right)=-\dfrac{c_+}{c_-}\,\phi\left(s_++\left(\dfrac{1}{c_+}-\dfrac{1}{c_-}\right)N\right)
\label{eq:4}
\ese
We substitute this into $\psi(s_+)+\phi(_+)=p(s_+)$ and get
\bsenn
-\dfrac{c_+}{c_-}\phi\left(s_++\left(\dfrac{1}{c_+}-\dfrac{1}{c_-}\right)N\right) +\phi(s_+)=p(s_+)
\esenn
Substituting $s_-$ back this gives
\bsenn
\phi(s_-) =  -\dfrac{c_-}{c_+}\,p\left(s_--\left(\dfrac{1}{c_+}-\dfrac{1}{c_-}\right)N\right)+
             \dfrac{c_-}{c_+}\, \phi\left(s_--\left(\dfrac{1}{c_+}-\dfrac{1}{c_-}\right)N\right)
\esenn
and that implies:
\bsenn
\phi(s_-)=-\sum_{k=1}^{\infty}\,\left(\dfrac{c_-}{c_+}\right)^k\,
 p\left(s_--\left(\dfrac{1}{c_+}-\dfrac{1}{c_-}\right)kN\right)
\esenn
Again this is a finite sum.

Summing $\psi(t-j/c_+)$ and $\phi(t-j/c_-)$ gives the general solution of the system (see Theorem \ref{theo:CV1})
\bsenn
z_j(t)  =  p\left(t-\dfrac{j}{c_+}\right) + \sum_{k=1}^\infty\,\left(\dfrac{c_-}{c_+}\right)^k\, \left[ p\left(t-\dfrac{j}{c_+}-\left(\dfrac{1}{c_+}-\dfrac{1}{c_-}\right)kN\right)- \right. \left. p\left(\left(t-\dfrac{j}{c_-}\right)-\left(\dfrac{1}{c_+}-\dfrac{1}{c_-}\right)kN\right) \right]
\esenn
Upon setting $j=N$, the terms telescope, and one obtains the solution
\bsenn
z_N(t)= \dfrac{c_+-c_-}{c_+} \sum_{k=0}^{\infty} \left(\dfrac{c_-}{c_+}\right)^k
p\left(t-\dfrac{N}{c_+}-\left(\dfrac{1}{c_+}-\dfrac{1}{c_-}\right)kN\right)
\esenn
Notice that according to Theorem \ref{theo:CV1},the functions $f_-$ and $f_+$ are not necessarily 
the same as the original pulse in the indicated time-interval (due to dispersion).  Hence we replace
$p$ by $f_k$.
$\blacksquare$

\vskip .0in\noindent
{\bf Remark:} In fact, from the proofs of Proposition 4.7 and Theorem 4.8 in \cite{CVH}, we can see 
that the high frequencies are quickly damped, and thus the $f_k$ in the above result are smoothed,
low-frequency versions of the original pulse.

The quantities used in the next Theorem are defined in Figure \ref{fig:schematic}.

\vskip .0in\noindent
\begin{theo} Let $S_N$ be as in Definition \ref{defn:boundary condns} and Equation \ref{eq:new-conds}.
The following holds independent of the boundary conditions. Fix a positive integer $k_0$.
As $N$ tends to infinity, then the orbit $z_N(t)$ of the last agent is characterized by
the following relations when $t\in[0,\frac{k_0N}{c_+}]$. Let $1\leq k\leq k_0$, then:
\bsenn
\begin{array}{ccc}
u_k &=& -\left(\dfrac{c_-}{c_+}\right)^kv_0\\
&&\\
\dfrac{T_k}{N} &=& \left(\dfrac{1}{c_+}-\dfrac{1}{c_-}\right)k\\
&&\\
\dfrac{A_k}{N} &=& -\left(\dfrac{c_-}{c_+}\right)^{k-1}\,\dfrac{v_0}{c_+}
\end{array}
\esenn
\label{theo:quantities}
\end{theo}

\vskip .0in\noindent
{\bf Proof:} Proposition \ref{prop:solution} gives us the solution of $P_N$, which is the
acceleration of the solution of $S_N$. This tells us that
at $t_k\equiv\dfrac{N}{c_+}+\left(\dfrac{1}{c_+}-\dfrac{1}{c_-}\right)kN$ velocities change by
\bsenn
\int_{t_k-\varepsilon}^{t_k+\varepsilon}\,z_N(t)\,dt =
\dfrac{c_+-c_-}{c_+}\left(\dfrac{c_-}{c_+}\right)^{k}v_0
\esenn
since by Principle 5, $f_k(t)$ is a unit pulse. Since the initial velocity is $u_0=-v_0$
(with respect to the leader), we obtain the following recursion for the velocities $u_k$ with
respect to the leader
\bsenn
u_{k+1}=u_k+\dfrac{c_+-c_-}{c_+}\left(\dfrac{c_-}{c_+}\right)^{k}\,v_0 \quad \logand \quad
u_0=-v_0
\esenn
and this gives the first result.

Integrating once more (and noting that $z_N(0)=0$), we see that the orbit $z_N(t)$ is given by a
piecewise affine function whose slope in the interval $(t_{k-1},t_k)$, with $t_{-1}=0$,
is given by $u_k$. From this we get the following recursion for the values of the local extrema
of $A_k$ of $z_N(t)$.
\bsenn
A_{k+1}=A_k+u_k\left(\dfrac{1}{c_+}-\dfrac{1}{c_-}\right)N \quad \logand \quad
A_1=-\dfrac{N}{c_+}v_0
\esenn
For the intercept times $T_k$ we set $z_N(T_k)=0$ and get
\bsenn
T_k= t_{k-1}+\tau_k \quad \where \quad  A_k+u_k\tau_k=0 \quad \logand \quad T_0=0
\esenn
$\blacksquare$

%\vskip .4in
\begin{centering}\section{Numerical Tests}
\label{chap:numerical}
\end{centering}
\setcounter{figure}{0} \setcounter{equation}{0}

To test these predictions (and the conjectures they are based upon), we ran 360
simulations and in each of these we measured 3 quantities.
For each $N\in\{100,200,400,800,1600,3200\}$ and each of the two boundary conditions
in Definition \ref{defn:boundary condns}, we took a grid of 30 parameter values: $\rho_{v,1}\in\{0,-0.1,-0.2,-0.3,-0.4,-0.5\}$ and $g_v\in\{-0.25,-0.50,-1,-2,-4\}$.
In accordance with Equation \ref{eq:new-conds}, this covers all the parameters still available
for variation in Equation \ref{eqn:oscillator array2}, since $\rho_{x,1}=-1/2$
and $|g_x|$ can be set equal to 1 by rescaling time (i.e. by introducing $\tau=\sqrt{|g_x|}t$).
We measured three quantities directly from numerical simulations that were done with
MATLAB's \texttt{ode45} algorithm with the relative and absolute tolerances set to $10^{-6}$:
the amplitude $A_1$ (see Figure \ref{fig:schematic}), the period $T$, which is the
(average) time elapsed between $T_k$ and $T_{k+2}$ (see Figure \ref{fig:schematic}), and the
attenuation $\alpha$, which is the ratio $A_{3}/A_1$ (see Figure
\ref{fig:schematic}). We then compared these with the predicted values obtained
in Theorem \ref{theo:quantities}. The relative errors (that is: $|$measured-predicted$|$/$|$predicted$|$)
are displayed in the figures below, one for each type of boundary conditions.
Note that each pair of points in the figures below represents measurements in each of the 30 grid points.
This is a log-log plot, so that the slope corresponds to the power of the decay.

\vskip .0in
The convergence is slow, as can be seen. It appears to be ${\cal O}(N^{-1/2})$ for the amplitude
and attenuation and ${\cal O}(N^{-1})$ for the period. It looks like it is not uniform
in the parameters nor in the boundary conditions.

\vskip .0in
%\begin{tabular}{|c|ll|ll|ll|}
%\hline
%REGULAR & Amplitude ($A_1$) &  & Period ($T$) &  & Attenuation ($\alpha$) & \\
%\hline
%Rel. Errors: & Average &Maximum &Average & Maximum &Average & Maximum\\
%\hline
%$N$= 100 & 0.04335 & 0.13355 & 0.00516 & 0.04502 & 0.07871 & 0.23406 \\
%\hline
%$N$= 200 & 0.03033 & 0.09466 & 0.00246 & 0.02390 & 0.05274 & 0.16194 \\
%\hline
%$N$= 400 & 0.02126 & 0.06701 & 0.00112 & 0.01195 & 0.03553 & 0.11182 \\
%\hline
%$N$= 800 & 0.01494 & 0.04741 & 0.00050 & 0.00542 & 0.02403 & 0.05763 \\
%\hline
%$N$=1600 & 0.01051 & 0.03354 & 0.00020 & 0.00209 & 0.01681 & 0.04021 \\
%\hline
%$N$=3200 & 0.00741 & 0.02372 & 0.00008 & 0.00064 & 0.01180 & 0.02818 \\
%\hline
%\end{tabular}
%
%\vskip .0in
%\begin{tabular}{|c|ll|ll|ll|}
%\hline
%VAR. MASS & Amplitude ($A_1$) &  & Period ($T$) &  & Attenuation ($\alpha$) & \\
%\hline
%Rel. Errors:  & Average &Maximum &Average & Maximum &Average & Maximum\\
%\hline
%$N$= 100 & 0.04136 & 0.12895 & 0.00599 & 0.04502 & 0.07860 & 0.23404 \\
%\hline
%$N$= 200 & 0.02932 & 0.09230 & 0.00282 & 0.02390 & 0.05271 & 0.16194 \\
%\hline
%$N$= 400 & 0.02076 & 0.06581 & 0.00128 & 0.01195 & 0.03552 & 0.11182 \\
%\hline
%$N$= 800 & 0.01468 & 0.04681 & 0.00057 & 0.00542 & 0.02403 & 0.05763 \\
%\hline
%$N$=1600 & 0.01038 & 0.03323 & 0.00023 & 0.00209 & 0.01681 & 0.04021 \\
%\hline
%$N$=3200 & 0.00734 & 0.02356 & 0.00009 & 0.00060 & 0.01179 & 0.02818 \\
%\hline
%\end{tabular}

\begin{figure}[pbth]
\center
\includegraphics[width=.8\linewidth]{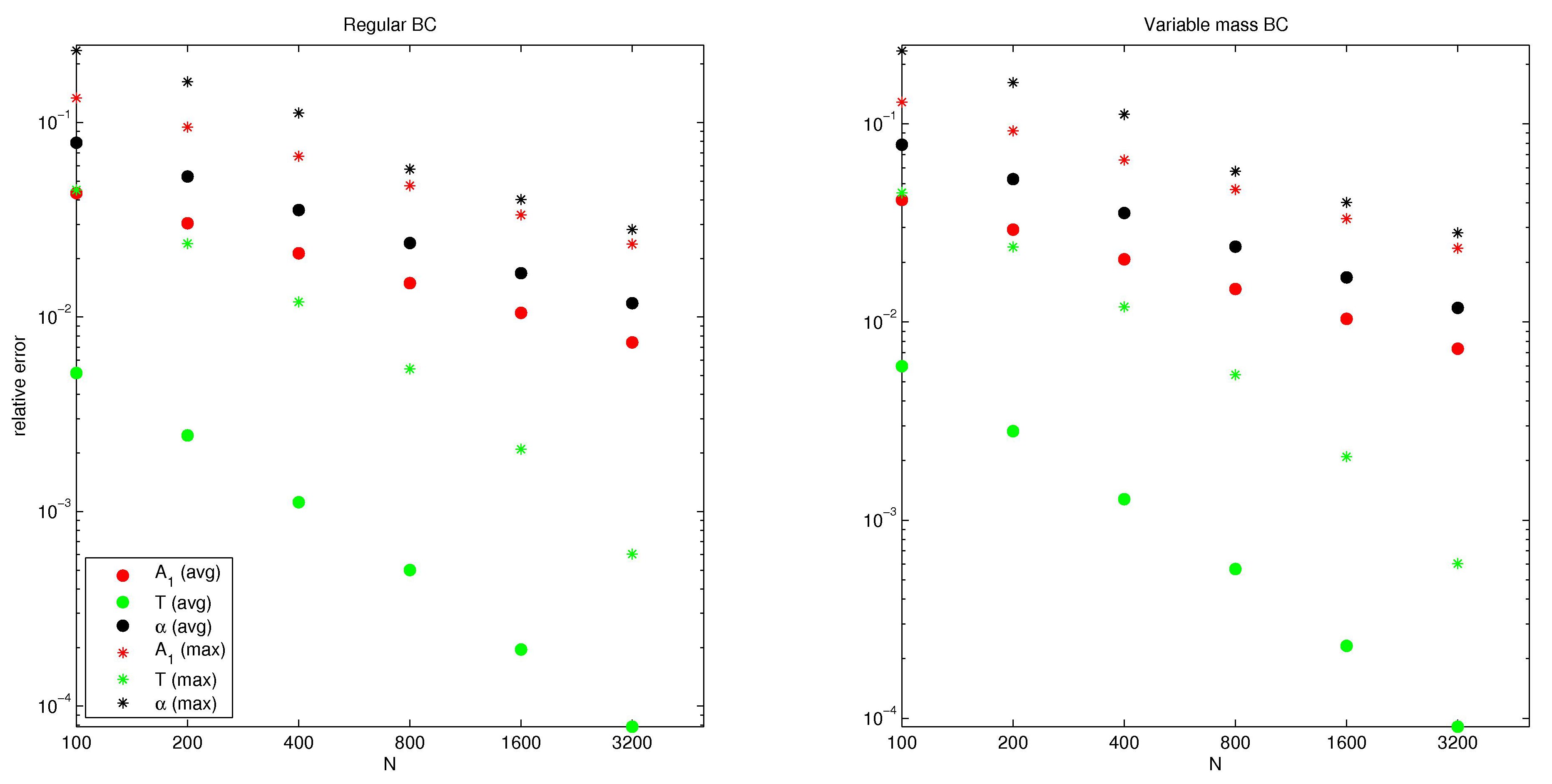} %for one-col
\caption{ \emph{Log-log plot of the relative errors (color online).
}}
\label{fig:errors}
\end{figure}

%\vskip .4in
\begin{centering}\section{Further Remarks}
 \label{chap:five}
 \end{centering}
\setcounter{figure}{0} \setcounter{equation}{0}

In this paper we have formulated and numerically checked a theory that gives an
explicit characterization of the transients for nearest neighbor systems of identical linear harmonic damped oscillators. The characterizations are asymptotic in the number of agents.
We established that transients will be smallest if the system behaves like the wave equation
and gave conditions on the parameters when that happens. In these cases we gave
precise quantitative characterizations of the transients in terms of the parameters of the problem.
These characterizations in this generality are new to the literature as is the fact that they
are independent of the boundary conditions.

Interestingly, it turns out that good strategies to obtain small transients,
involve non-symmetric interactions. Other authors (\cite{BMH}, \cite{HMHS}, and \cite{LFJ}))
have also noticed this phenomenon, but they did so in the context of slightly different
models that are not strictly decentralized. An example of that is the model obtained from
Equation \ref{eqn:oscillator array2} by replacing
$ g_v \left\{\rho_{v,-1}\dot z_{k-1} + \rho_{v,0}\dot z_k +\rho_{v,1} \dot z_{k+1}\right\}$
with $g_v\dot z_k$.

\vskip .0in\noindent
Now that we know the transient as function of the parameters, it is reasonable to look
for an optimal choice of parameters. There are many quantities that could be optimized.
To illustrate our point here, we simply look at an index that involves the squares of the amplitudes,
and by using Theorem \ref{theo:quantities} we get
\bsenn
I_E\equiv\dfrac{1}{N^2 v_0^2}\,\sum_{k=1}^\infty A_k^2 =  \dfrac{1}{c_+^2-c_-^2}
\esenn
Recall that according to that theorem the amplitude of the
oscillations actually \emph{increase} exponentially if $|c_-|>|c_+|$ (given in
Equation \ref{eq:signal-velocities}). It is easy to check that this happens
whenever $\rho_{v,-1}\leq-1/2$. A straightforward calculation shows the following.

\begin{figure}[pbth]
\center
\includegraphics[width=.4\linewidth]{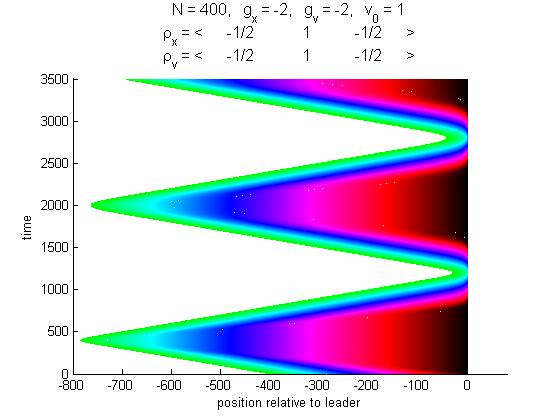}
\includegraphics[width=.4\linewidth]{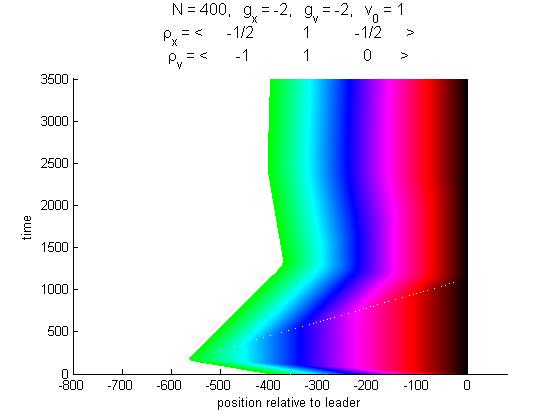}
\caption{ \emph{Simulations for $N=400$, $g_x=g_v=-2$ and $\rho_{x,1}=-1/2$, regular boundary conditions.
In the first picture $\rho_{v,1}=-1/2$, in the second $\rho_{v,1}=0$.
Each color represent the orbit of one of the 400 individual agents.
Horizontal is position relative to the leader, vertical is time.
}}
\label{fig:transient2}
\end{figure}

%\noindent
Note that $g_x=0$ gives a good coefficient ($I_E=1$) but does not spatially align the agents.
As an example choose $g_x=g_v=-2$, $\rho_{x,1}=-1/2$,
and let $\rho_{v,1}$ take the values $-1/2$ and 0. We take $N=400$ and the leader leaves with velocity $v_0=1$. The simulations are exhibited in Figure \ref{fig:transient2}.
The table below gives predictions for the various quantities. It is clear that the
asymmetric choice $\rho_{v,1}=0$ far outperforms the symmetric choice $\rho_{v,1}=-1/2$.

%\vskip .0in\noindent
\begin{center}
\begin{tabular}{|c|c|c|c|c|c|c|}
\hline
 & $c_+$ & $c_-$ & $ A_1$ & $T$ & $\alpha$ & $I_E$ \\
\hline
$\rho_{v,1}=-1/2$ & 1 & -1 & 400 & 1600 & 1 & $\infty$\\
\hline
$\rho_{v,1}=0$ & $1+\sqrt{2}$ & $1-\sqrt{2}$ & 166 & 2262 & 0.029 & 0.177\\
\hline
\end{tabular}
\end{center}

%\vskip 0.4 in

\vspace{\fill}
\end{document}